\documentclass[reqno]{amsart}     
\usepackage{pstricks,pst-node,pst-coil,pst-plot}   

\theoremstyle{plain}      
 
\newtheorem*{thmhr}{Theorem (Hebey, Robert; [HR04])}

\newtheorem{thm}{Theorem}[section]     
\newtheorem{theorem}[thm]{Theorem}

\theoremstyle{remark}      
 
\newtheorem{remark}[thm]{Remark}

\theoremstyle{definition}

\def\al{{\alpha}}

\def\om{{\omega}}

\def\si{{\sigma}}

\def\ep{{\epsilon}}

\def\Th{{\Theta}}         
         
\def\phi{{\varphi}}

\let\pa\partial

\DeclareMathAlphabet{\doba}{U}{msb}{m}{n}

\gdef\mN{\doba{N}}

\gdef\mR{\doba{R}}         
\gdef\mS{\doba{S}}

\def\ric{{\mathop{\rm Ric}}}

\def\div{{\mathop{\rm div}}}     

\def\grad{{\mathop{\rm grad}}}

\def\d2dt{\frac{d^2}{dt^2}}
\def\grad{{\mathop{\rm grad}}}

\def\eref#1{{\rm (\ref{#1})}}

\let\<\langle 
\let\>\rangle

\long\def\ignorethis#1{}

\newdimen\templaenge

\def\Atbox#1#2{\setbox0\hbox{$\displaystyle #1$}\templaenge=\textwidth\advance\templaenge by -\wd0%
\setbox1\hbox{$#2$}\advance\templaenge by -\wd1%
$$#1\hbox{\kern\templaenge$#2$\hss}$$\par\bigbreak}

\newtheorem*{subcaseI.1}{Subcase I.1}
\newtheorem*{subcaseI.2}{Subcase I.2}

\newtheorem*{subcaseII.1}{Subcase II.1}
\newtheorem*{subsubcaseII.1.1}{Subsubcase II.1.1}
\newtheorem*{subsubcaseII.1.2}{Subsubcase II.1.2}
\newtheorem*{subcaseII.2}{Subcase II.2}

\newcommand{\definedas}{\mathrel{\raise.095ex\hbox{\rm :}\mkern-5.2mu=}}

\begin{document}     


\title{Positive mass theorem for the Paneitz-Branson operator}


\author{Emmanuel Humbert} 
\address{Institut \'Elie Cartan, BP 239 \\ 
Universit\'e de Nancy 1 \\
54506 Vandoeuvre-l\`es-Nancy Cedex \\ 
France}
\email{humbert@iecn.u-nancy.fr}

\author{Simon Raulot} 
\address{Institut de Math\'ematiques\\
Universit\'e de Neuch\^atel\\
Rue Emile-Argand 11 \\2007 Neuch\^atel\\ Suisse}
\email{simon.raulot@unine.ch}


\begin{abstract}
We prove that under suitable assumptions, the constant term in the Green
function of the Paneitz-Branson operator on a compact Riemannian manifold $(M,g)$
is positive unless $(M,g)$ is conformally diffeomophic to the standard
sphere. The proof is inspired by the positive mass theorem on spin
manifolds by Ammann-Humbert \cite{ammann.humbert:03}.
\end{abstract}


%
%



\maketitle     




\section{Introduction}


Let $(M,g)$ be a compact Riemannian manifold of dimension $n \geq 4$. We denote by $Q_g$ the  $Q$-curvature for the metric $g$ defined by 
$$ Q_g \definedas \frac{n^2 - 4}{8n(n-1)^2}S_g^2 - \frac{2}{(n-2)^2}|E_g|^2 + \frac{1}{2(n-1)} \Delta_g S_g,$$ 


\noindent where $\Delta_g = - \div_g \nabla$ is the Laplace-Beltrami operator, $S_g$ stands for the scalar curvature of $g$, $|E_g|$ denotes the
$g$-norm of the Einstein tensor $E_g\definedas \ric_g- \frac{S_g}{n} g$ and $\ric_g$ is the Ricci curvature of $g$. The Paneitz-Branson 
operator introduced for $n = 4$ by Paneitz  in \cite{paneitz:83} and whose definition was generalized in dimension greater than $5$ by Branson \cite{branson:87}, is defined for all $u \in C^\infty(M)$ by 
$$P_g u \definedas \Delta_g^2 u - \div_g\left(A_g du \right) + \frac{n-4}{2} Q_gu$$

\noindent where 
$$A_g \definedas \frac{(n-2)^2+4}{2(n-1)(n-2)}S_g g - \frac{4}{n-2}\ric_g.$$

\noindent This operator is closely related to the problem of prescribing $Q$-curvature in a conformal class as well as the
Yamabe operator  (see \eref{yamabe_equation} below) is related to the problem of prescribing the scalar curvature in a conformal class. It is a
conformally covariant operator  in the sense that if $g' = e^{2f} g$ is conformal to $g$, then for all $v \in C^\infty(M)$, 
$$P_{g'}( e^{- \frac{n-4}{2}f} v) = e^{-\frac{n+4}{2}f}P_g (v).$$ 

\noindent In particular, if $n \geq 5$, and if we set $u =e^{\frac{n-4}{2}f}$ so that $g'  = u^{\frac{4}{n-4}} g$, we
get for all $v \in C^\infty(M)$
\begin{eqnarray} \label{conformal_cov1}
P_{g'}(u^{-1} v) = u^{-\frac{n+4}{n-4}} P_g(v).
\end{eqnarray}
 
\noindent From now on, we make the following assumptions: \\

\noindent {\bf (a)} $g$ is conformally flat;

\noindent {\bf (b)} $n \geq 5$;

\noindent {\bf (c)} the Yamabe invariant is positive (see for instance
\cite{aubin:98} or \cite{hebey:97}) i.e. $g$ is conformal to a metric
$g'$ for which the scalar curvature is positive.

\noindent {\bf (d)} the operator $P_g$ is positive. \\

\noindent  Under Assumptions ${\bf (a)}$ to ${\bf (d)}$, it is well known that the Green's function
$G_g$ of $P_g$ exists, is unique and smooth on $M \setminus\{ p\}$. By the conformal convariance of the Paneitz-Branson 
operator, if $g' = u^{\frac{4}{n-4}} g $ is conformal to $g$, then 
$$G_{g'} (x,y) = \frac{G_g(x,y)}{u(x)u(y)}.$$ 

\noindent Now, let $p \in M$. By $(1)$, up to a conformal change of metric, we can assume \\

\noindent {\bf (a')} $g$ is flat around $p$. \\

\noindent Then, it is known that we have the following expansion when $x$ is close to $p$,  
\begin{eqnarray} \label{exp1}
G_g(x,p) = \frac{1}{2(n-2)(n-4)\om_{n-1} d_g(x,p)^{n-4} } + A +\alpha_p(x)
\end{eqnarray}

\noindent where $\om_{n-1}$ stands for the volume of the $(n-1)$-dimensional sphere, $A \in \mR$, $\al_p$ is a smooth function defined
around $p$ and satisfying $\al_p(p) = 0$. By analogy to the case of the conformal Laplacian (see again \cite{aubin:98,hebey:97}), 
the number $A$ is called the {\it mass of the Paneitz-Branson operator}. If $g'= u^{\frac{4}{n-4}}g$ is another
metric conformal to $g$ and flat around $p$, then the mass $A'$ corresponding to the metric $g'$ is given by  
$$A' = \frac{A}{u(p)^2}.$$

\noindent Hence, the mass $A$ depends on the choice of the metric in the
conformal class, but not its sign. This is the reason why in  the statement of Theorem
\ref{main} below, we do not need to assume ${\bf (a')}$.  \\

\noindent  We also  make the following assumption \\

\noindent {\bf (e)} $G_g > 0$ on $M \setminus \{p \}$. \\

\noindent  For interesting results concerning Assumptions ${\bf (d)}$ and
${\bf (e)}$, the reader may refer to Grunau-Robert 
\cite{grunau.robert:07}. \\

\noindent  The main result of the paper is the following:

\begin{theorem} \label{main}
Under assumptions ${\bf (a)}$ to ${\bf (e)}$, the mass $A$ satisfies 
$$A\geq 0$$
with equality if and only if $(M,g)$ is conformally diffeomorphic to the sphere. 
\end{theorem}
 
\noindent Theorem \ref{main} has been already proven with the additional assumption that the Poincar\'e exponent is small enough (see
\cite{qing.raske:06a, qing.raske:06b}). In this case, Qing and Raske proved also the positivity of the Green's function of $G_g$. 

\noindent Our proof is inspired from the positive mass theorem on spin manifolds by Ammann-Humbert in \cite{ammann.humbert:03} 
(see also Raulot \cite{raulot:07}). The difficulty  here is to overcome the fact that on non-spin manifolds, there is no equivalent of the 
Schr\"odinger-Lichnerowicz Formula.\\

\noindent Hebey and Robert proved  the nice
following result which is an analogue for geometric equations of order 4  of a hard
problem concerning the Yamabe Equation:     
\begin{thmhr}
Let $(M,g)$ be a conformally flat compact manifold of dimension $n \geq 5$. Assume $g$ has a positive Yamabe invariant, that $P_g$ is positive as well as its Green function and that the mass of $P_g$ is positive. Then, the geometric equation
$$P_gu = u^{\frac{n+4}{n-4}}$$
is compact.   
\end{thmhr}

\noindent In particular, together with Theorem \ref{main}, we get rid of the positivity of the mass.\\

\noindent {\bf Acknoledgements} We  want to thank Emmanuel Hebey and Fr\'ed\'eric Robert who gave us  many helpful informations and
references on the subject.


\section{Proof of Theorem \ref{main}} \label{proof}


\noindent In the whole proof, we can work with Assumption ${\bf (a')}$
which does not restrict the generality as explained above. To avoid complicated formulas, we set 
$$H(x) = 2(n-2)(n-4)\om_{n-1} G_g(x,p).$$

\noindent By Relation \eref{exp1}, $H$ satisfies the following expansion near $p$ 
\begin{eqnarray} \label{exp2}
H(x) = \frac{1}{d_g(x,p)^{n-4}} + B + \al(x)
\end{eqnarray}

\noindent where $B = 2(n-2)(n-4)\om_{n-1} A$ and where  $\al =  2(n-2)(n-4)\om_{n-1}\al_p$ is smooth around $p$ and satisfies
$\al(p)=0$. Theorem \ref{main} is equivalent to show that $B \geq 0$ with equality if and only if $(M,g)$ is conformally 
diffeomorphic to the standard sphere. 
 
\noindent For any metric $g$, let 
$$L_g \definedas \frac{4(n-1)}{n-2} \Delta_g + S_g$$

\noindent be the Yamabe operator. We recall some well known facts about $L_g$. The reader may refer to \cite{aubin:98,hebey:97} for further informations. First, as well as $P_g$, $L_g$ is conformally covariant. If $g'= u^\frac{4}{n-2} g$ is conformal to $g$ then 
\begin{eqnarray} \label{conformal_cov2}	
L_{g'} (u^{-1}\,\cdot\,) = u^{-\frac{n+2}{n-2}} L_g(\,\cdot\,)
\end{eqnarray}

\noindent It follows that the scalar curvatures $S_g$ and $S_{g'}$ are related by the following equation
\begin{eqnarray}\label{yamabe_equation}
L_g u = S_ {g'} u^{\frac{n+2}{n-2}}.
\end{eqnarray}
    
\noindent By Assumptions ${\bf (a')}$ and ${\bf (b)}$, the Green's function $\Lambda_g$ of $L_g$ exists, is unique, smooth 
and positive on $M \setminus \{ p \}$. Setting $\Gamma(x) = 4 (n-1) \om_{n-1}\Lambda_g(x,p)$ to simplify formulas, we have 
when $x$ is close to $p$ 
\begin{eqnarray} \label{exp3}
\Gamma(x) = \frac{1}{d_g(x,p)^{n-2}} + C + \beta(x)
\end{eqnarray}

\noindent where by $C \in \mR$, $\beta$ is a smooth function defined around $p$ and satisfies $\beta(0)= 0$. 
We define a new metric  $g' \definedas \Gamma^{\frac{4}{n-2}} g$ conformal to $g$ on $M_0 \definedas M \setminus \{ p\}$. 
Then, by \eref{yamabe_equation} 
\begin{eqnarray} \label{S=0}
S_{g'} = \Gamma^{-\frac{n+2}{n-2}} L_g(\Gamma) \equiv 0
\end{eqnarray}
 
\noindent on $M_0$. We set $H'= \Gamma^{-\frac{n-4}{n-2}} H$. By conformal covariance of the Paneitz-Branson operator 
\eref{conformal_cov1} and since $P_g H = 0$ on $M_0$, we have $P_{g'} H' \equiv 0 \hbox{ on } M_0$.  Define for all $\ep >0$ 
small enough, $M_\ep \definedas M \setminus B^g(p,\ep)$ where $B^g(p,\ep)$ stands for the ball of center $p$ and radius $\ep$ 
with respect to the metric $g$. We have 
\begin{eqnarray} \label{intP1}
\int_{M_\ep} P_{g'} H' dv_{g'} = 0.
\end{eqnarray}

\noindent By Relation \eref{S=0} and from the definition of $P_g$ we have 
$$P_{g'}H' = \Delta_{g'}^2 H' - \div_{g'}\left( \frac{4}{n-2} \ric_{g'} dH' \right) -\frac{n-4}{(n-2)^2}|E_{g'}|^2H'.$$ 

\noindent Set $S_{\ep} \definedas \partial M_{\ep}= \partial B^g(p,\ep)$ be the $(n-1)$-dimensional sphere of center $p$ and radius $\ep$.
We let $ds_{g'}$ (resp. $ds_g$) be the volume element induced by $g'$ (resp. g) on $S_\ep$. Integrating by part the above relation, we obtain 
\begin{eqnarray} \label{intP2}
\int_{M_\ep} P_{g'} H' dv_{g'}  = -{\bf I} + \frac{4}{n-2}\,{\bf II} - \frac{1}{2} \int_{M_\ep} |E_{g'}|^2 H' dv_{g'} 
\end{eqnarray}

\noindent where
\[ \left\{ \begin{array}{ccc}
{\bf I} & =&  \int_{S_\ep} \partial_{\nu'} \Delta_{g'} H' ds_{g'} \\
{\bf II} & = &  \int_{S_\ep} \ric_{g'}(\grad^{g'}H ', \nu') ds_{g'}. 
\end{array}
\right. \]

\noindent Here, $\nu'$ denotes the unit outer normal vector on $S_\ep = \partial M_\ep$ with respect to the metric $g'$.


\subsection{Computation of ${\bf I}$}


First, we notice that the scalar curvatures $S_g$ and $S_{g'}$ vanish on $S_\ep$. For $g$, this comes from Assumption ${\bf (a')}$ 
and for $g'$, this follows from \eref{S=0}. Consequently, using Formula \eref{conformal_cov2} and 
\begin{eqnarray*}
\Delta_{g'} H' & =&  \frac{n-2}{4(n-1)} L_{g'} H'\\
& = & \frac{n-2}{4(n-1)} \Gamma^{-\frac{n+2}{n-2}} L_{g} \left( \Gamma H' \right)\\
& = &  \frac{n-2}{4(n-1)} \Gamma^{-\frac{n+2}{n-2}} L_{g} \left( \Gamma^{\frac{2}{n-2}} H \right)
\end{eqnarray*}

\noindent We obtain 
\begin{eqnarray} \label{laplacians}
\Delta_{g'} H'  =  \Gamma^{-\frac{n+2}{n-2}} \Delta_g\left( \Gamma^{\frac{2}{n-2}} H \right).
\end{eqnarray}

\noindent We set $r \definedas d_g(x,p)$. From Formulas \eref{exp2} and \eref{exp3}, we have: 
\begin{eqnarray*} 
\Gamma^{\frac{2}{n-2}} H  = \left(\frac{1}{r^{n-2}} + C+ \beta(x) \right)^{\frac{2}{n-2} } \left( \frac{1}{r^{n-4}} + B+ \al(x)\right).
\end{eqnarray*}

\noindent Then, using Taylor formula at $p$,
$$\Gamma^{\frac{2}{n-2}} H   = r^{2-n} + B r^{-2} + O(r^{-1}) $$

\noindent where in the whole proof, $O(r^m)$ denotes a smooth function defined in a neighborhood of $p$ and which satisfies  
\begin{eqnarray*} 
|\nabla_g^k O(r^m)|_g \leq C_k r^{m-k}
\end{eqnarray*}

\noindent for all $k\in \mN$. Since $g$ is flat around $p$, we have for radially symmetric functions $f$,
\begin{eqnarray} \label{polar_coor}
\Delta_g f(r) = -f''(r) -\frac{n-1}{r}f'(r). 
\end{eqnarray}

\noindent Hence, this gives that near $p$,
$$\Delta_g \Gamma^{\frac{2}{n-2}} H = 2(n-4) B r^{-4} + O(r^{-3})$$

\noindent and hence by \eref{laplacians} and \eref{exp3}
$$\Delta_{g'} H' = \Gamma^{- \frac{n+2}{n-2}} \Delta_{g} H = 2(n-4)B r^{n-2} + O(r^{n-1}). $$

\noindent We then obtain 
\begin{eqnarray} \label{ddr}
\frac{\partial}{\partial r} \left( \Delta_{g'} H'\right) = 2(n-2)(n-4)B r^{n-3} + O(r^{n-2}).   
\end{eqnarray}
 
\noindent On $S_\ep$, $r \equiv \ep$. In addition, 
\begin{eqnarray} \label{nu}
\nu' =- \Gamma^{-\frac{2}{n-2}} \frac{\partial}{\partial r} =  - (\ep^2 + o(\ep^2)) \frac{\partial}{\partial r}
\end{eqnarray}

\noindent and 
\begin{eqnarray} \label{volume_element}
ds_{g'} = \Gamma^{2\frac{n-1}{n-2}} ds_g = \Gamma^{2\frac{n-1}{n-2}} \ep^{n-1}ds = \left(\ep^{1-n} + o(\ep^{1-n}) \right)ds.
\end{eqnarray}

\noindent where $ds$ stands for the standard volume element on the unit $(n-1)$-sphere. 
By Formulas \eref{ddr}, \eref{nu} and \eref{volume_element}, we obtain 
\begin{eqnarray} \label{I}
{\bf I} & =& -2(n-2)(n-4) \om_{n-1} B + o(1) 
\end{eqnarray}


\subsection{Computation of ${\bf II}$}


If $g' = e^{2f} g$ is conformal to $g$, then the following formula holds (see \cite{hebey:97} p. $240$ or \cite{aubin:98}):
\begin{eqnarray} \label{ricg'}
\ric_{g'} = \ric_g - (n-2) \nabla^2 f + (n-2) \nabla f \otimes \nabla f +
\left(\Delta_g f  -(n-2) |\nabla f|^2_g\right) g.   
\end{eqnarray}

\noindent In this context, $f= \frac{2}{n-2} \log(\Gamma)$. By \eref{exp3}, we have near $p$ 
\begin{eqnarray} \label{f}  
f &=  & \frac{2}{n-2}  \log\left( \frac{1}{r^{n-2}} + O(1)  \right)
\nonumber \\
& = &  -2 \log(r) +  O(r^{n-2}). 
\end{eqnarray}

\noindent Let $(r, \Th_1,\cdots,\Th_{n-1})$ be polar coordinates on $\mR^n$. The Christoffel symbols $\Gamma^r_{r,\Th_i}$ of 
the Euclidean metric in these coordinates identically vanish. This implies that for any radially symmetric function $h$, the 
mixed terms $\nabla^2_{r \Th_i}h$ are zero. Since $g$ is flat near $p$, we deduce that 
\begin{eqnarray*}
\nabla^2 f & = & \frac{2}{r^2} dr^2+ b+
\bar{O}(r^{n-4}) 
\end{eqnarray*}

\noindent where, as in what follows, we denote by $\bar{O}(r^m)$ a 2-form whose norm with respect to $g$ is $O(r^m)$ and where 
$b$ is a $2$-form such that 
\begin{eqnarray} \label{b}
b\left(\,\cdot\,,\frac{\pa}{\pa r} \right)\equiv 0.
\end{eqnarray}  

\noindent Using \eref{polar_coor}, one also computes that 
\begin{eqnarray*}
\nabla f  \otimes \nabla f & =& \frac{4}{r^2} dr^2 +\bar{O}(r^{n-4})\\
\Delta_g f & = & \frac{2(n-2)}{r^2} + O(r^{n-4}) \\
| \nabla f|^2_g& =&  \frac{4}{r^2} + O(r^{n-4}).
\end{eqnarray*}

\noindent Since $g$ is flat near $p$, $\ric_g$ vanishes and $g = dr^2 + r^2 \sigma^{n-1}$ where $\si^{n-1}$ stands for the usual 
metric on the standard sphere $\mS^{n-1}$. We deduce from these computations that 
\begin{eqnarray*}
\ric_{g'} & = &  -(n-2)b  -\frac{2(n-2)}{r^2} dr^2 
 +   \frac{4(n-2)}{r^2} dr^2 +\\
& & \left(   \frac{2(n-2)}{r^2} - 
 \frac{4(n-2)}{r^2} +O(r^{n-4}) \right)(dr^2 + r^2
\si^{n-1} )+ \bar{O}(r^{n-4})\\
& = & -(n-2)b -2(n-2) \sigma^{n-1} + \bar{O}(r^{n-4}).
\end{eqnarray*}

\noindent We get from \eref{exp2}, \eref{exp3} and the definition of $H'$ that on $S_\ep$    
\begin{eqnarray} \label{grad}
\grad^{g'}H' &  =& \Gamma^{-\frac{4}{n-2}} \grad^g\left( 1 + O(r^{n-4})
\right) 
\nonumber \\
& = &  O(r^{n-1}) \frac{\partial}{\partial r} + v 
\end{eqnarray}

\noindent is a vector field such that  $\ric_{g'}(v, \nu')=0$. Observe that by \eref{nu} and \eref{b},
we have $\sigma^{n-1}(\,\cdot\,,\nu')=0$  and $b(\,\cdot\,,\nu')=0$ on $S_\ep$. In addition, the estimates \eref{nu}, \eref{b} then 
imply that on $S_\ep$
\begin{eqnarray*} 
 \ric_{g'}(\grad^{g'} H',\nu')& =&  \bar{O}(r^{n-4})(\grad^{g'} H', \nu')  \\
& = &   O(\ep^{2n-3}).
\end{eqnarray*}

\noindent Relation \eref{volume_element} then leads to
\begin{eqnarray} \label{II}
{\bf II} = O(\ep^{n-2})= o(1). 
\end{eqnarray}


\subsection{Conclusion}


Using  \eref{intP1}, \eref{intP2}, \eref{I}, \eref{II} and passing to the limit $\ep \to 0$, we obtain that 
\begin{eqnarray} \label{finalngeq6} 
0 = 2(n-2)(n-4) \om_{n-1} B - \frac{1}{2} 
\int_{M\setminus\{p\}} |E_{g'}|^2 H' dv_{g'}.
\end{eqnarray}

\noindent Assumption ${\bf (e)}$ implies that $H' > 0$ and hence $B \geq 0$. This proves first part of Theorem \ref{main}.

\noindent Now, assume that $B =0$. Then $E_{g'} \equiv 0$ on $M \setminus \{p\}$. This implies that  $(M \setminus \{p\}, g')$ 
is Einstein and scalar flat hence Ricci flat. Since in addition the Weyl
curvature is zero, $(M\setminus \{p\}, g')$ turns to  
be flat (see \cite{hebey:97} p. $123$) . It is known that $(M\setminus
\{p\}, g')$ is asymptotically flat and that its mass satisfies $m(g') =
c_n  C$  where $c_n >0$ (see e.g. Lee-Parker \cite{lee.parker:87}). Since $g'$ is flat,
$m(g')=0$ so is $C$ and by a positive mass Theorem by Schoen-Yau \cite{schoen.yau:88}, $(M,g)$ is conformally diffeomorphic to $(S^n,g)$.

\begin{remark} 
\noindent It is clear from the proof that Assumption ${\bf (a)}$ can  be weakened and replaced by 

\noindent {\bf (a)} $g$ is locally flat around a point $p$ and the standard
Positive Mass Theorem is valid on $M$ (i.e. with the notations of
Section \ref{proof}, $C \geq 0$ with equality if and only if $(M,g)$ is
conformally diffeomorphic to $\mS^n$). In particular, by
\cite{schoen.yau:79} and \cite{ammann.humbert:03}, this assumption holds
if $n\in \{ 5,6,7 \}$ or if $M$ is spin.

\end{remark}


\end{document}